\documentclass[article,12pt]{elsarticle}
\usepackage{lineno,hyperref}
\modulolinenumbers[5]
\usepackage{lineno}
\modulolinenumbers[5]
\usepackage{amssymb}
\usepackage{color,soul}
\usepackage{grffile}
\usepackage{booktabs}
\usepackage{times}
\usepackage{subcaption}
\usepackage{amsmath}
\usepackage{amsthm}
\usepackage{ulem}
\usepackage[titletoc]{appendix}
\usepackage{titletoc}
\usepackage{mathtools}
\usepackage{doi}
\usepackage{url}
\usepackage{tikz}
\usepackage{graphicx}
\usepackage{graphics,color}
\usepackage{makeidx}
\usetikzlibrary{calc}
\usepackage{multirow}
\usetikzlibrary{patterns}
\usepackage[margin=1.0in]{geometry}
\usepackage[outdir=./]{epstopdf}
\usepackage{euscript}
\usepackage{natbib}
\usepackage{wrapfig}
\usepackage{epsfig}

\theoremstyle{definition}

\usetikzlibrary{matrix, shapes, arrows, positioning, chains}

\begin{document}

\begin{frontmatter}
\title{A Single Equation Explains Go-or-Grow Dynamics in Cyclic Hypoxia}
\author[ad]{Gopinath Sadhu\corref{mycorrespondingauthor}}
\cortext[mycorrespondingauthor]{Corresponding author}
\ead{gopinaths@iisc.ac.in}
\author[f]{Philip K Maini}
\ead{philip.maini@maths.ox.ac.uk}
\author[ad]{Mohit Kumar Jolly}		

\ead{mkjolly@iisc.ac.in}
\address[ad]{Department of Bioengineering, Indian Institute of Science, Bangalore, Karnataka, India}
\address[f]{ Wolfson Centre for Mathematical Biology, Mathematical Institute, University of Oxford, Oxford, United Kingdom}			
\begin{abstract}
We propose a minimal mathematical framework to describe the go-or-grow dynamics of tumor cells comprising two phenotypically distinct populations. One population is migratory and undergoes linear diffusion, while the other proliferates in an oxygen-dependent manner. The local oxygen concentration governs transitions between these phenotypes. We then ask whether these two coupled phenotype-specific equations can be reduced to a single mixed-phenotype equation under cyclic hypoxia. We establish a connection between the minimal go-or-grow model with distinct phenotypic populations and a reduced model describing a single-cell population with oxygen-dependent diffusion and proliferation in the fast-phenotypic-switching regime. This theoretical reduction is validated through numerical simulations.
\end{abstract}
\begin{keyword}
\texttt Go-or-grow, cyclic hypoxia, phenotype switching
\end{keyword}
		
	\end{frontmatter}
	
	
\section{Introduction}
The go-or-grow dichotomy is a well-established mathematical modeling hypothesis whereby cell populations either migrate or proliferate, but not both simultaneously. This modeling framework has been used to describe many biological and medical conditions, mainly the tumor invasion process \cite{2024crossley,2012_andreas,2025go_grow_hillen}. Under the go-or-grow hypothesis, cells can switch between proliferative and migratory states via phenotypic switching, as observed in experiments \cite{vittadello2020examining}. This phenotypic switching may depend on many factors, such as cell density and local oxygen levels \cite{2024crossley, 2024Bekar, sadhu2026phenotype}. Local oxygen levels in the tumor microenvironment may fluctuate between normoxic (this environment promotes migratory to proliferative phenotype switching) and hypoxic levels (this environment encourages proliferative to migratory phenotypic transition), which is termed as cyclic hypoxia \cite{sadhu2026phenotype,2019cyclic_MKJ}. 

In this work, based on this set of observations, we propose a minimal go-or-grow model under cyclic hypoxia using a system of reaction-diffusion equations. The proposed model consists of two equations for proliferative and migratory cells: migratory cells move via linear diffusion, proliferative cells reproduce based on local oxygen concentration, and phenotypic switching is a function of local oxygen level. The goal of this work is to capture go-or-grow dynamics through a single mixed phenotype equation model instead of the proposed two equation model for heterogeneous cell phenotypes. The model has three timescales: phenotypic switching, cell proliferation and migration, and oxygen fluctuations. 
In the fast phenotypic switching regime, we show that the model involving distinct phenotypic populations can be reduced to a model that consists of a single mixed-phenotype population, in which the diffusion coefficient and proliferation rate are linked to the phenotypic switching rates and proliferation rate of the two-population model. Our numerical results confirm that the single mixed-phenotype cell model behavior is the same as that of the proposed go-or-grow model with heterogeneous populations for a specific choice of biologically meaningful phenotypic switching functions. In the rapid oxygen fluctuation regime, where the timescale of oxygen fluctuations is very much shorter than those of phenotype switching and proliferation/movement, we observe, however, that the single-mixed-phenotype equation fails to reproduce the two-equation go-or-grow model behaviour. 
\section{A two population go-or-grow model}
Let $m(\mathbf{x},t)$ and $p(\mathbf{x},t)$ denote the densities of migratory and proliferative cells, respectively, at time $t$ and position $\mathbf{x}=(x,y,z)$ in cartesian coordinate space and define the total cell density by
\begin{equation}
n(\mathbf{x},t) = m(\mathbf{x},t) + p(\mathbf{x},t).
\end{equation}
A minimal go-or-grow model, where phenotypic switching and proliferation depend on oxygen concentration ($c(t)$), is given by
\begin{align}
     &\frac{\partial p}{\partial t}=\mu(c) p\left(1-\frac{n}{K}\right)-\lambda_1(c) p +\lambda_2 (c) m, \label{eq:proliferative}\\
      &\frac{\partial m}{\partial t}=D_m\nabla^2m
    +\lambda_1(c) p -\lambda_2 (c) m \label{eq:migratory},
\end{align}
where $\nabla^2=\frac{\partial^2}{\partial x^2}+\frac{\partial^2}{\partial y^2}+\frac{\partial^2}{\partial z^2}$. In most in vitro experiments done so far to recapitulate cyclic hypoxia, spatial heterogeneity has not been the focus, instead temporal fluctuations of oxygen levels have been imposed on cells \cite{sadhu2026EMT,sadhu2026phenotype,celora2024characterising}. Hence, $ c (t)$ fluctuates over time between oxygen-rich and oxygen-poor levels. 
Here, the first term on the right-hand side of Eq. \eqref{eq:proliferative} describes the proliferative cells that proliferate logistically with rate $\mu(c)$, which is an increasing function of $c$. Proliferative cells adopt a migratory phenotype with a rate $\lambda_1(c)$, which is decreasing with increasing $c$. The $\lambda_2(c)$ denotes migratory cells switching to proliferation and is an increasing function of $c$. 
Here, and unless stated otherwise, we keep $\mu(c)$, $\lambda_1(c)$, and $\lambda_2(c)$ general, assuming them to be strictly non-negative functions of the oxygen concentration $c$.
\subsection{A single-equation formulation for the coupled phenotypes model}\noindent
A compact way to express the go-or-grow dynamics under cyclic hypoxia is by the single equation 
\begin{equation}
\frac{\partial n}{\partial t}=\nabla \cdot \big( D(c) \nabla n \big)+r(c)\, n\left( 1 - \frac{n}{K} \right),
\end{equation}
The go-or-grow trade-off is encoded through the condition
\begin{equation}
D'(c)< 0~\text{and}~ r'(c) >0
\end{equation}
which ensures that increased proliferation is accompanied by reduced motility with increasing $c$.\\ \noindent\\
Our goal is to derive expressions for $D(c)$ and $r(c)$ in terms of $\mu(c)$, $\lambda_1(c)$, and $\lambda_2(c)$.

\subsection{Connecting $D(c)$ and $r(c)$ with $\mu(c)$, $\lambda_1(c)$, and $\lambda_2(c)$}
Experimental studies reported that the time-scale of switching for both proliferative to migratory and migratory to proliferative phenotypes can be faster than that of cellular-level dynamics, such as cell proliferation and cell migration \cite{celia2018hysteresis}. We define the parameter $\varepsilon$ to be the ratio between the phenotypic switching time scale and the population dynamics time scale, namely
\begin{equation}
\varepsilon=\frac{\text{phenotypic switching time scale}}{\text{population dynamics time scale}}.
\end{equation}
The condition $0<\varepsilon \ll 1$ implies that phenotypic adaptation occurs much faster
than migration or proliferation.
Fast phenotypic switching is modeled for Eqs. \eqref{eq:proliferative}-\eqref{eq:migratory}  by assuming \cite{2024Bekar}
\begin{equation}
\lambda_1(c) = \frac{1}{\varepsilon}\,\tilde{\lambda}_1(c),
\qquad
\lambda_2(c) = \frac{1}{\varepsilon}\,\tilde{\lambda}_2(c),
\qquad
0 < \varepsilon \ll 1,
\end{equation}
where $\tilde{\lambda}_1(c), \tilde{\lambda}_2(c)$ are $\mathcal{O}(1)$.\\ \noindent\\
Hence, the typical fast-switching system takes the form
\begin{align}
&\frac{\partial p}{\partial t}=\mu(c)\, p \left( 1 - \frac{n}{K} \right)+ \frac{1}{\varepsilon}\tilde{\lambda}_2(c)\, m
- \frac{1}{\varepsilon}\tilde{\lambda}_1(c)\, p,\label{eq:fast_proli}\\
&\frac{\partial m}{\partial t}
=\nabla \cdot \big( D_m \nabla m \big)- \frac{1}{\varepsilon}\tilde{\lambda}_2(c)\, m + \frac{1}{\varepsilon}\tilde{\lambda}_1(c)\, p\label{eq:fast_migra}.
\end{align}
In the singular limit $\varepsilon \to 0$, the leading-order balance yields the relation
\begin{equation}
\tilde{\lambda}_2(c)\, m
=
\tilde{\lambda}_1(c)\, p,
\end{equation}
 and the two cell densities $p$ and $m$ can be related to the total cell density ($n$) via
\begin{equation}
m = \theta(c)\, n,
\qquad
p = \big(1 - \theta(c)\big)\, n,~\text{where}
~
0<\theta(c)
=
\frac{\lambda_1(c)}
{\lambda_1(c) + \lambda_2(c)}<1.
\end{equation}
Substituting these expressions into the total density equation (sum of Eqs. \eqref{eq:fast_proli} and \eqref{eq:fast_migra})  leads to the reduced
single-equation model
\begin{equation}\label{eq:mixed_single}
\frac{\partial n}{\partial t}=\nabla \cdot \big( D(c) \nabla n \big)+r(c)\, n \left( 1 - \frac{n}{K} \right),
\end{equation}
with effective oxygen-dependent coefficients
\begin{equation}
D(c) = \theta(c)\, D_m,
\qquad
r(c) = \big(1 - \theta(c)\big)\, \mu (c).
\end{equation}
\subsection{The conditions $D'(c)<0$ and $r'(c)>0$ are satisfied for the mixed phenotype single equation}\noindent
Let $A=(\lambda_1(c)+\lambda_2(c))$, then
\begin{align}
    \theta'(c)=\frac{\lambda_1'(\lambda_1+\lambda_2)-\lambda_1(\lambda_1'+\lambda_2')}{A^2}\nonumber
    =\frac{\lambda_1'\lambda_2-\lambda_1\lambda_2'}{A^2}<0~\text{as}~\lambda_1>0,\lambda_2>0, \lambda_1'<0~\text{and}~\lambda_2'>0.
\end{align}
Hence, $D'(c)=D_m\theta'(c)<0$ and $r'(c)=(1-\theta(c))\mu'(c)-\theta'(c)\mu(c)>0$ as $\mu>0$ and $\mu'>0$.
\section{Numerical results}\noindent
We take the functional forms of the proliferation rate and phenotype switching rates to be
\begin{equation}\label{eq:rates_assume}
 \mu(c)=\mu_p\Psi(c),~   \lambda_1(c)=\lambda_{pm}(1-\Psi(c))~\text{and}~ \lambda_2(c)=\lambda_{mp}\Psi(c)~\text{with}~\Psi(c)=\frac{c^k}{c_H^k+c^k},
\end{equation}
where $\mu_p$, $\lambda_{pm}$, and $\lambda_{mp}$ represent the rates of cell proliferation, the transition from a proliferative to a migratory phenotype, and the transition from a migratory to a proliferative phenotype, respectively. Additionally, $c_H$ denotes the hypoxic threshold value of oxygen concentration, and $k$ refers to the Hill coefficient.
Here $\mu'(c)>0,$ $\lambda'_2(c)>0$ and $\lambda'_1(c)<0$. Hence, $\theta(c)=\frac{\lambda_{pm}c_H^k}{\lambda_{pm}c_H^k+\lambda_{mp}c^k}$ yields $D(c)=D_m\frac{\lambda_{pm}c_H^k}{\lambda_{pm}c_H^k+\lambda_{mp}c^k}$ and $r(c)=\frac{\lambda_{mp}c^k}{\lambda_{pm}c_H^k+\lambda_{mp}c^k}\big(\frac{\mu_p c^k}{c_H^k+c^k}\big)$. To obtain numerical results, we employed the forward time central difference in space (FTCS) scheme within a finite difference framework in a one-dimensional Cartesian domain $[0, L]$ under no-flux boundary conditions for Eq. \eqref{eq:migratory}. For the heterogeneous model (Eqs. \eqref{eq:proliferative} and \eqref{eq:migratory}), we assumed that the tumor initially grows from proliferative cells only. Hence, the initial conditions are considered as $p(x,0)=0.2 \mathcal{H}(x-2)$ and $m(x,0)=0$  for $x \in [0, L]$ with $L=10$. For the mixed phenotype model (Eq. \ref{eq:mixed_single}), we use homogeneous Neumann condition at the boundary of the domain $[0, L]$ and initial condition as $n(x,0)=0.2 \mathcal{H}(x-2)$ for $x\in [0, L]$. During cyclic hypoxia, \begin{equation}
    c(t)=\begin{cases}
        c_1>c_H, \text{in}~ (2q-2)\frac{{\mathrm{T_{oxy}}}}{2}<t\leq (2q-1)\frac{{\mathrm{T_{oxy}}}}{2}\\
         c_2<c_H, \text{in}~ (2q-1)\frac{{\mathrm{T_{oxy}}}}{2}<t\leq (2q)\frac{{\mathrm{T_{oxy}}}}{2}, 
    \end{cases}
\end{equation}
where $T_{\mathrm{oxy}}$ is the oxygen fluctuation period and $q=1,\dots,s$ with $sT_{\mathrm{oxy}}$ being the total time of interest.
\begin{figure}[h!]
	\centering
	\includegraphics[scale=0.62]{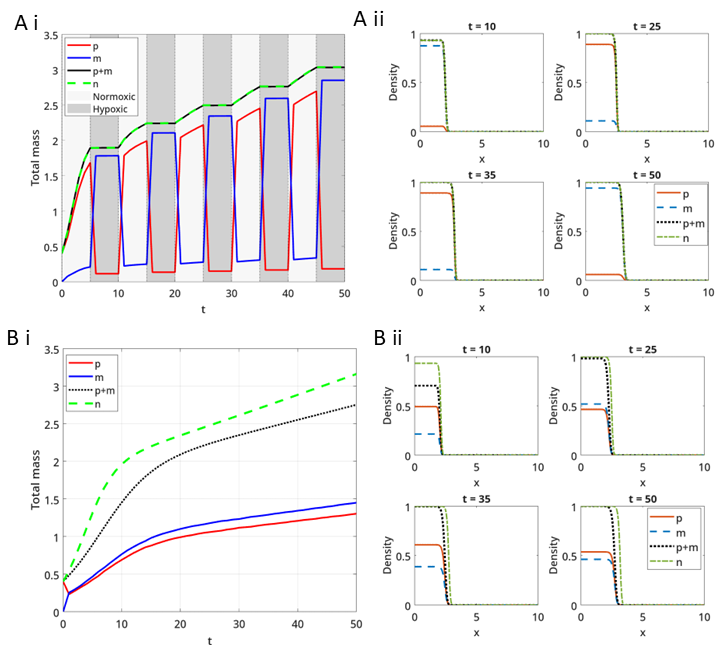}
	\caption{Time evolution of the proliferative ($p$), migratory cell ($m$) and total (sum of proliferative and migratory cells) and mixed phenotype cell density ($n$) in (A) fast phenotype regimes and (B) rapid oxygen oscillation regimes. Here, $D_m=0.01$, $\lambda_{pm}=10$,$\lambda_{mp}=20$, $\mu_p=1$, $L=10, K=1, c_H=0.5$ and the hill coefficient $k=3$. The period of cyclic hypoxia ($T_{\mathrm{oxy}}$) is 10 in (A) and 0.001 in (B).}
	\label{fig:figure1}
\end{figure}	
In this problem, we have three timescales, namely $T_{\mathrm{switch}}=\max\{\frac{1}{\lambda_{mp}},\frac{1}{\lambda_{pm}}\}$ denotes the phenotypic adaptation time scale, $T_{\mathrm{pop}}=\max\{\frac{1}{\mu_p}, \frac{L^2}{D_m}\}$ the population dynamics time scale associated with diffusion and proliferation, and oxygen fluctuation timescale ($T_{\mathrm{oxy}}$). In the fast phenotype switching regime i.e., $T_{\mathrm{switch}}\ll T_{\mathrm{oxy}}<T_{\mathrm{pop}}$, we observed that the total density of the distinct heterogeneous population model follows the same time course evolution as for the mixed phenotype heterogeneous model (Figure \ref{fig:figure1}A). Under this condition, phenotypic switching equilibrates before significant migration or growth occurs. However, in the rapid oxygen oscillation regime i.e., $T_{\mathrm{oxy}} \ll T_{\mathrm{switch}}<T_{\mathrm{pop}}$, cells are unable to adapt to instantaneous oxygen variations. As a result, rapid periodic oxygen fluctuations do not induce go-or-grow alternation in this regime (Figure \ref{fig:figure1}B).
\section{Model analysis for limiting cases}\noindent
We simplify the switching function $\Psi(c)$ to be (in the limit $k\rightarrow\infty$) $\Psi(c) = \mathcal{H}(c-c_H)$, where $\mathcal{H(.)}$ is the Heaviside function:
\begin{equation*}
    \mathcal{H}(c-c_H)=\begin{cases}
       & 1; c>c_H \\
       & 0; \text{otherwise}
    \end{cases}
\end{equation*}
Hence, the rates for proliferation, the transition from proliferative to migratory, and the transition from migratory to proliferative for Eqs. \eqref{eq:proliferative}-\eqref{eq:migratory} given in Eqs. \eqref{eq:rates_assume}, reduce to
\begin{equation}\label{eq:reduced_rates_assume}
 \mu(c)=\mu_p\mathcal{H}(c-c_H),~   \lambda_1(c)=\lambda_{pm}\mathcal{H}(c_H-c)~\text{and}~ \lambda_2(c)=\lambda_{mp}\mathcal{H}(c-c_H).
\end{equation}
We now consider the case where oxygen levels are constant and, as before, we assume that the tumor grows from proliferative cells only, so, initially, $p=p_0$ and $m=0$, hence $n=p_0$.\\ \noindent
\textbf{For oxygen-rich conditions, i.e., $c>c_H$:} In this case, proliferative cells cannot transition into migratory cells so, from the initial conditions, the migratory cell density remains at zero. As a result, the two system model reduces to:
\begin{align}\label{eq:normoxic_reduced}
     \frac{\partial p}{\partial t}=\mu_p p\left(1-\frac{p}{K}\right).
\end{align}
Here $\theta(c)=0$, so, $D(c)=0$ and $r(c)=\mu_p$. Hence, the single mixed equation \eqref{eq:mixed_single} becomes
\begin{equation}\label{eq:mixed_reduced}
     \frac{\partial n}{\partial t}=\mu_p n\left(1-\frac{n}{K}\right).
\end{equation}
Clearly, Eq. \eqref{eq:normoxic_reduced} and Eq. \eqref{eq:mixed_reduced} are identical.\\
\noindent
\textbf{For oxygen-deprived conditions, i.e., $c<c_H$:} The two system model becomes
\begin{align}
     &\frac{\partial p}{\partial t}=-\lambda_{pm}p, \label{eq:proli_hypoxia}\\
     &\frac{\partial m}{\partial t}=D_m\nabla^2m+\lambda_{pm}p. \label{eq:migration_hypo}
\end{align}
By adding Eqs. \eqref{eq:proli_hypoxia} and \eqref{eq:migration_hypo}, we have 
\begin{equation}
    \frac{\partial n}{\partial t}=D_m\nabla^2m
\end{equation}
However, under this oxygen-deprived condition, all proliferating cells (p) switch to migratory cells (m) rapidly, as we are in the fast phenotypic switching regime. Hence, $m$ is equivalent to $p+m$ i.e., $n$.\\ \noindent
In this case, $\theta(c)=1$, so, $D(c)=D_m$ and $r(c)=0$. Therefore, the single mixed equation \eqref{eq:mixed_single} becomes
\begin{equation}
    \frac{\partial n}{\partial t}=D_m\nabla^2n,
\end{equation}
which is exactly same as the sum of Eqs. \eqref{eq:proli_hypoxia} and \eqref{eq:migration_hypo} in fast phenotype switching regime.
\section{Discussion and future direction}
In this work, we analyzed a minimal go-or-grow model, in which cells can either proliferate via logistic growth or move via linear diffusion, with oxygen-dependent proliferation rate and phenotypic switching rates. We established a relationship between the model of two phenotypic cell types (``specialist'') and a model for a single mixed equation in the fast-phenotype-switching regime under cyclic hypoxia. This relation suggests that the diffusion coefficient and proliferation rate for the mixed model are determined by oxygen-dependent switching rates between phenotypically distinct cell population models. Our simulation results show that in the fast phenotype-switching regime, the mixed-phenotype model closely resembles a two-phenotype population model. We have also shown that if the oxygen oscillation timescale is very rapid (for example near blood vessels), then the mixed-phenotype model fails to accurately capture the go-or-grow hypothesis.

In this work, we have considered the oxygen levels to be spatially homogeneous but, in reality, they will vary in space, as well, as time \cite{murphy2023growth}. Future work will aim to determine if a single mixed-phenotype equation can capture the dynamics of a two distinct phenotype model in this more complex situation. Furthermore, in this work, we have only
considered linear diffusion of migratory cells, however, more complex models have been proposed in which migration is affected by total cell density and the extracellular matrix \cite{2024crossley}, giving, potentially, a more realistic picture of tumor invasion.
\section*{Data availability}
The code for numerical simulations are available at GitHub.
\section*{Acknowledgment}
GS and MKJ are supported by Param Hansa Philanthropies. PKM would like to thank the Royal Society Yusuf Hamied Visiting Fellowship scheme which funded his visit to IISc Bangalore. For the purpose of open access, the author has applied a CC BY public copyright licence to any author accepted manuscript arising from this submission.
\bibliography{AML}

@article{2024crossley,
  title={Phenotypic switching mechanisms determine the structure of cell migration into extracellular matrix under the ‘go-or-grow’hypothesis},
  author={Crossley, Rebecca M and Painter, Kevin J and Lorenzi, Tommaso and Maini, Philip K and Baker, Ruth E},
  journal={Mathematical Biosciences},
  pages={109240},
  year={2024},
  publisher={Elsevier}
}

@article{celora2024characterising,
	title={Characterising cancer cell responses to cyclic hypoxia using mathematical modelling},
	author={Celora, Giulia L and Nixson, Ruby and Pitt-Francis, Joe M and Maini, Philip K and Byrne, Helen M},
	journal={Bull. Math. Bio.},
	volume={86},
	number={12},
	pages={145},
	year={2024},
	publisher={Springer}
}

@article{murphy2023growth,
  title={Growth and adaptation mechanisms of tumour spheroids with time-dependent oxygen availability},
  author={Murphy, Ryan J and Gunasingh, Gency and Haass, Nikolas K and Simpson, Matthew J},
  journal={PLOS Computational Biology},
  volume={19},
  number={1},
  pages={e1010833},
  year={2023},
  publisher={Public Library of Science San Francisco, CA USA}
}

@article{sadhu2026phenotype,
  title={A phenotype-structured PDE framework for investigating the role of hypoxic memory on tumor invasion under cyclic hypoxia},
  author={Sadhu, Gopinath and Jain, Paras and George, Jason Thomas and Jolly, Mohit Kumar},
  journal={Bulletin of Mathematical Biology},
  volume={88},
  number={2},
  pages={23},
  year={2026},
  publisher={Springer}
}

@article{sadhu2026EMT,
  title={The impact of oxygen heterogeneity on epithelial-mesenchymal transitions: a numerical study},
  author={Sadhu, Gopinath and Byrne, Helen M and Dalal, DC},
  journal={Journal of Mathematical Biology},
  volume={92},
  number={1},
  pages={15},
  year={2026},
  publisher={Springer}
}

@article{2019cyclic_MKJ,
  title={Acute vs. chronic vs. cyclic hypoxia: their differential dynamics, molecular mechanisms, and effects on tumor progression},
  author={Saxena, Kritika and Jolly, Mohit Kumar},
  journal={Biomolecules},
  volume={9},
  number={8},
  pages={339},
  year={2019},
  publisher={MDPI}
}

@article{2024Bekar,
  title={Travelling waves in a minimal go-or-grow model of cell invasion},
  author={Falc{\'o}, Carles and Crossley, Rebecca M and Baker, Ruth E},
  journal={Applied Mathematics Letters},
  volume={158},
  pages={109209},
  year={2024},
  publisher={Elsevier}
}

@article{2025go_grow_hillen,
  title={Go-or-grow models in biology: a monster on a leash},
  author={Thiessen, Ryan and Conte, Martina and Stepien, Tracy L and Hillen, Thomas},
  journal={Journal of Mathematical Biology},
  volume={91},
  number={5},
  pages={58},
  year={2025},
  publisher={Springer}
}

@article{celia2018hysteresis,
  title={Hysteresis control of epithelial-mesenchymal transition dynamics conveys a distinct program with enhanced metastatic ability},
  author={Celi{\`a}-Terrassa, Toni and Bastian, Caleb and Liu, Daniel D and Ell, Brian and Aiello, Nicole M and Wei, Yong and Zamalloa, Jose and Blanco, Andres M and Hang, Xiang and Kunisky, Dmitriy and others},
  journal={Nature communications},
  volume={9},
  number={1},
  pages={5005},
  year={2018},
  publisher={Nature Publishing Group UK London}
}

@article{vittadello2020examining,
  title={Examining go-or-grow using fluorescent cell-cycle indicators and cell-cycle-inhibiting drugs},
  author={Vittadello, Sean T and McCue, Scott W and Gunasingh, Gency and Haass, Nikolas K and Simpson, Matthew J},
  journal={Biophysical Journal},
  volume={118},
  number={6},
  pages={1243--1247},
  year={2020},
  publisher={Elsevier}
}

@article{2012_andreas,
    author = {Hatzikirou, H. and Basanta, D. and Simon, M. and Schaller, K. and Deutsch, A.},
    title = {‘Go or Grow’: the key to the emergence of invasion in tumour progression?},
    journal = {Mathematical Medicine and Biology: A Journal of the IMA},
    volume = {29},
    number = {1},
    pages = {49-65},
    year = {2010},
    month = {07}
}
\bibliographystyle{unsrt}
\end{document}